\documentclass[10pt,twoside]{article}
\usepackage{latexsym}
\usepackage{amssymb}

\makeatletter

\@addtoreset{equation}{section} \makeatother
\newtheorem{theorem}{Theorem}[section]
\newtheorem{lemma}[theorem]{Lemma}
\newtheorem{proposition}[theorem]{Proposition}

\newtheorem{definition}[theorem]{Definition}
\newtheorem{remark}[theorem]{Remark}

\newcommand{\g}{\langle\cdot,\cdot\rangle }
\newcommand{\m}{{\cal M}}
\newcommand{\mo}{{\cal M}_0}
\newcommand{\J}{{\cal J}}

\newcommand{\es}{{\cal S}}

\newcommand{\dimo}{\noindent{\it Proof.\ }}
\newcommand{\cvd}{\hfill \rule{0.5em}{0.5em}}
\newcommand{\N}{{\mathbb N}}

\newcommand{\R}{{\mathbb R}}
\newcommand{\LL}{{\mathbb L}}
\newcommand{\NN}{\Omega^1_K(p,q)}
\begin{document}
\title{{\bf\Large Global hyperbolicity and Palais--Smale condition for
action functionals in stationary spacetimes}}
\author{{\bf\large A.M. Candela\footnote{{\sl Corresponding author}.
Supported by M.I.U.R. (research funds ex 40\% and 60\%).} ,
J.L. Flores\footnote{Supported by MEC Grants MTM-2004-06262 and
RyC-2004-382.} ,
M. S\'anchez\footnote{Supported by MEC Grant MTM-2004-04934-C04-01.}}\\
\\
{\it\small $^*$Dipartimento di Matematica, Universit\`a degli Studi di Bari,}\\
{\it\small Via E. Orabona 4, 70125 Bari, Italy}\\
{\it\small candela@dm.uniba.it}\\
{\it\small $^\dagger$Departamento de \'Algebra, Geometr\'{\i}a y Topolog\'{\i}a,}\\
{\it \small Facultad de Ciencias, Universidad de M\'alaga,}\\
{\it\small Campus Teatinos, 29071 M\'alaga, Spain}\\
{\it\small floresj@agt.cie.uma.es}\\
{\it\small $^\ddagger$Departamento de Geometr\'{\i}a y Topolog\'{\i}a,}\\
{\it\small Facultad de Ciencias, Universidad de Granada,}\\
{\it\small Avenida Fuentenueva s/n, 18071 Granada, Spain}\\
{\it\small sanchezm@ugr.es}}

\date{}

\maketitle

\centerline{\sl Dedicated to the memory of Professor Jerzy
Konderak}

\smallskip

\begin{center}
{\bf\small Abstract} \vspace{3mm} \hspace{.05in}\parbox{4.5in}
{{\small In order to apply variational methods to the action
functional for geodesics of a stationary spacetime, some
hypotheses, useful to obtain classical Palais--Smale condition,
are commonly used: pseudo--coercivity, bounds on certain
coefficients of the metric, etc. We prove that these technical
assumptions admit a natural interpretation for the conformal
structure (causality) of the manifold. As a consequence, any
stationary spacetime with a complete timelike Killing
vector field and a complete Cauchy hypersurface (thus,
globally hyperbolic), is proved to be geodesically connected.}}
\end{center}



\section{Introduction}

In the last years,  an intensive research on the problem of
geodesic connectedness in stationary spacetimes (i.e., the
question whether any two points in a Lorentzian manifold
admitting a  timelike Killing vector field, can be
joined by a geodesic) has been carried out. Even though there are
geometric and physical reasons --- no analog to Hopf--Rinow
theorem exists for a Lorentzian manifold, stationary spacetimes
include typical physical spacetimes, as Kerr's or Schwarzschild's
--- the main interest comes from the analytical viewpoint. In
fact, given a Lorentzian manifold  $(\m,\g_{L})$,
geodesics connecting two fixed points $p,q\in \m$ are, among
the $C^1$ curves connecting them, the critical points of the
(energy) action functional
\begin{equation}
\label{funzionaleazione} f(z) =\ \frac{1}{2}\int_0^1\langle\dot
z,\dot z\rangle_L\ ds ,
\end{equation}
which becomes strongly indefinite in the Lorentzian setting.
Moreover, authors have mainly followed an extrinsic approach based
on the assumption that the spacetime is {\em standard stationary},
i.e. $(\m,\g_{L})$ splits globally as ${\cal M} =\mo\times\R$,
with $({\cal M}_0,\langle\cdot,\cdot\rangle)$ a finite dimensional
connected Riemannian manifold, and metric $\g_{L}$ written as
\begin{equation}\label{met}
\langle \cdot,\cdot\rangle_L= \langle\cdot,\cdot\rangle + 2
\langle\delta (x),\cdot\rangle dt - \beta(x)dt^2
\end{equation}
for each $T_z{\cal M} \equiv T_x\mo \times \R$, $z = (x,t) \in
{\cal M}$, where $\delta$ and $\beta$ are a smooth vector field
and a smooth strictly positive scalar field on $\mo$,
respectively. In this case, the lack of boundedness of $f$ can be
overcome by means of a suitable variational principle, stated in
\cite{BF, BFG, giama1}, which shows that looking for critical
curves of action functional $f$ connecting $p=(x_p,t_p)$ to
$q=(x_q,t_q)$ becomes equivalent to the study of critical points
for a new functional $\J$ on the Riemannian part, namely
\begin{equation}\label{newfunction}
\begin{array}{rcl}\displaystyle
\J(x) &=&\displaystyle \frac 12\int_0^1\langle\dot x,\dot
x\rangle\ ds\
 +\ \frac 12\int_0^1\frac{\langle \delta(x),\dot x\rangle^2}{\beta(x)}\ ds\\
&&\; \displaystyle\ -\ \frac 12\ \left(\int_0^1 \frac{\langle
\delta(x),\dot x\rangle}{\beta(x)}\ ds - \Delta_t\right)^2 \
\left(\int_0^1 \frac1{\beta(x)}\ ds\right)^{-1}
\end{array}
\end{equation}
($\Delta_t = t_q-t_p$), which is defined on a suitable set of
``spatial'' curves joining $x_p$ to $x_q$ in $\mo$ (for more
details, see Proposition \ref{newazione}) and may also be bounded
from below.

Since then, such a functional has been widely studied. Considering
only the case $\g$ complete (without boundary), the known main
results can be summarized as follows:
\begin{itemize}
\item[1.] Benci, Fortunato and Giannoni \cite{BF, BFG} studied the
geodesic connectedness in a standard static spacetime
($\delta\equiv 0$) and introduced functional $\J$ for this case.
Giannoni and Masiello \cite{giama1} extended this study to the
standard stationary case. From these results (see also
\cite[Theorem 3.4.3]{Ma}), the existence of critical points of
$\J$ is ensured when $\beta$ and
$|\delta(x)|^2=\langle\delta(x),\delta(x)\rangle$ have an upper
bound and some $\epsilon\in \R$ exists so that
\begin{equation}\label{below}
0<\epsilon \leq \beta(x) \qquad \hbox{for all $x \in \mo$}.
\end{equation}
\item[2.] Pisani \cite{Pi} used a different approach based on the
direct study of action functional $f$. He obtained that under
assumption (\ref{below}) a sublinear growth for $\beta$ and
$\delta$ suffices for the existence of critical points, i.e., it
is enough to assume that some $\alpha<1$ exists so that
\begin{equation}\label{sublin2}
\beta(x), | \delta(x)| \le \mu d^{\alpha}(x,\bar x) + k \quad
\mbox{for all $x \in \mo$,}
\end{equation}
where $d(\cdot,\cdot)$ is the canonical distance associated to
$\g$, $\bar x\in \m_0$ is fixed, and $\mu \ge 0$, $k \in \R$ (see
\cite[Theorem 1.2]{Pi} and also \cite{CS} for a multiplicity
result). 
\item[3.] Remarkably, Giannoni and Piccione \cite{GP}
studied the existence of critical points for action $f$ from a
more intrinsic viewpoint. In principle, they assume only the
existence of a complete timelike Killing vector field $K$. Then,
for each $p,q\in \m$, they introduce a natural space of curves
$C^1_K(p,q)$ associated to $K$, and consider the restriction of
$f$ to this space. Then, they define a notion of {\em
pseudo--coercivity} for $f$, and show that, under this condition,
$p$ and $q$ can be joined by a geodesic. However, there are two
important limitations on this result: (a) pseudo--coercivity
implies global hyperbolicity and, thus, the spacetime must be
isometric to a standard stationary one (see next section), and (b)
more unpleasantly, pseudo--coercivity is a very technical
analytical condition. So, in order to give a more concrete result,
they fix a Cauchy temporal function (which becomes equivalent to
choose a splitting (\ref{met})), and reprove Pisani's result (at
least when $\beta = -g(K,K)$ is bounded). \item[4.] R. Bartolo and
the authors \cite{BCFS} have applied very accurate estimates (some
of them coming from \cite{CFS}) to functional $\J$
 in the standard static case, i.e., when the
functional in (\ref{newfunction}) is simplified by $\delta\equiv
0$. As a consequence, the exponent $\alpha=2$ in (\ref{sublin2})
is shown to be enough and optimal for the existence of critical
points of $\J$. Recently, Bartolo, Candela and Flores \cite{BCF}
have extended this result to the stationary case, showing that it
is sufficient to assume
\begin{eqnarray}\label{quad}
0\ <\ \beta(x) &\le& \mu_1 d^{2}(x,\bar x) + k_1 \quad \mbox{for
all $x \in \mo$,}\\ \label{lineare} |\delta (x)| &\le& \mu_2
d(x,\bar x) + k_2\ \quad \mbox{for all $x \in \mo$}
\end{eqnarray}
(with $\bar x \in \m_0$ fixed, and $\mu_1, \mu_2 \ge 0$, $k_1,k_2
\in \R$). We must emphasize that these hypotheses correspond to
the rough bounds in order to ensure the global hyperbolicity of
the spacetime (see Appendix A.1).
\end{itemize}

The study of the standard stationary case hides an important fact:
the same spacetime can split as (\ref{met}) in very different ways
(with very different $\beta, \delta$) because just one such
splitting is not intrinsic to the spacetime. As a simple and
extreme example, Minkowski spacetime can be written as (\ref{met})
either with an arbitrary growth of $| \delta|$ or with an
incomplete $\g$ (see Appendix A.2). More deeply, the bounds for
$\beta$, $| \delta|$ do not have a geometric meaning on $\m$,
except as sufficient (but neither necessary nor intrinsic)
conditions for global hyperbolicity.

Motivated by this type of objections, here we focus on Giannoni
and Piccione's approach. The main limitation of their results is
that pseudo--coercivity condition is analytical and very
technical. In fact, it can be regarded as a tidy and neat version
of Palais--Smale condition for the stationary ambient. But now the
question is how to translate this technical condition in terms of
the (Lorentzian) geometry of the manifold.

The aim of this paper is to answer this question by showing that,
essentially, the geometrical meaning of pseudo--coercivity is {\em
global hyperbolicity with a complete Cauchy hypersurface}. In
fact, we prove the following result (which extends all the
previous ones), discussing carefully all the hypotheses:
\begin{theorem}\label{globhyp}
Let $(\m,\g_{L})$ be a stationary spacetime
 with a complete
timelike Killing vector field $K$. If $\m$ is globally hyperbolic
with a complete (smooth, spacelike) Cauchy hypersurface $\es$,
then it is geodesically connected.
\end{theorem}
Even more, well--known standard arguments in previous references
(based on Ljusternik--Schnirelman category) allow one to prove
also some multiplicity results when $\m$ is not contractible in
itself. Concretely: (i) each two points $p,q\in \m$ can be joined
by a sequence of spacelike geodesics $(z_n)_n$ with diverging
$(f(z_n))_n$, (ii) given $p\in \m$ and an integral curve $\gamma$
of $K$, the number of (future-directed) timelike geodesics which
connect $p$ and $\gamma(s)$, diverges when $s\rightarrow +\infty$.

This paper is organized as follows. In Section \ref{section2} we
introduce the necessary Lorentzian tools, emphasizing the
interplay between stationarity and global hyperbolicity. In
Section \ref{sec3} Giannoni--Piccione's intrinsic approach is
revisited, and the relevant aspects for our problem are stressed.
As the hypotheses of Theorem \ref{globhyp} will imply the
existence of a splitting as in (\ref{met}), in Section \ref{sec4}
we explain how the functional approach is simplified when  one
chooses one such splitting (but the results will be independent of
the chosen one). The proof of Theorem \ref{globhyp} is carried out
in Section \ref{sec5}. Previous discussion translates it into
Theorem \ref{coerciveness}, and the crucial step for its proof is
Proposition \ref{p-s}. Essentially, this proposition shows that,
because of the existence of a complete Cauchy hypersurface, one
must check Palais--Smale condition only for sequences of curves
with bounded Riemannian norm. In the Appendix A we give exhaustive
examples and discussions which show the accuracy of Theorem
\ref{globhyp} and explain the meaning of the involved hypotheses.
In general, we try to minimize the technicalities in the interplay
between Causality Theory and Variational Methods (see Remark
\ref{lipschitz2}). Nevertheless, one of these technicalities,
which regards the relation between continuous and $H^1$ causal
curves is interesting in its own right, and is studied in Appendix
B.


\section{Tools in Lorentzian Geometry}\label{section2}

In this section we briefly recall some basic notions in Lorentzian
Geometry which will be used along the paper (for more details on
Lorentzian manifolds, see \cite{BEE, HE, ON, Pe, SW}).

By a {\em Lorentzian manifold} $(\m,\g_L)$ we mean a
smooth\footnote{As a simplification, {\em smooth} will mean
$C^\infty$ as usual. But this hypothesis can be relaxed. In fact,
one can assume that smooth means only $C^4$ for the spacetime and,
consequently, $C^3$ for the elements which depend on first
derivatives, as hypersurfaces. This allows one to apply Nash
Embedding Theorem (in the spirit of most previous references on
this topic) to some hypersurfaces, and no new problem will appear
in our case, because the existence of smooth {\em Cauchy}
hypersurfaces has been proved in \cite{BeS1} (see discussion
below).} (connected) finite dimensional manifold equipped with a
semi--Riemannian metric of index 1 on each tangent space $T_z\m$,
$z\in\m$. A tangent vector $\zeta \in T_z{\cal M}$ is called {\em
timelike} (respectively {\em lightlike}; {\em spacelike}; {\em
causal}) if $\langle \zeta, \zeta \rangle_L <0 $ (respectively
$\langle \zeta, \zeta \rangle_L = 0$ and $\zeta \neq 0$; $\langle
\zeta, \zeta \rangle_L >0 $ or $\zeta =0$; $\zeta$ is either
timelike or lightlike). In what follows the Lorentzian manifold
$(\m,\g_L)$ will be also a {\em spacetime}, that is, $(\m, \g_L)$
is connected and time--orientable, with a prescribed
time--orientation (a continuous choice of a causal cone at each
$p\in \m$, which is called the {\em future} cone, in opposition to
the non--chosen one or {\em past} cone).

A $C^1$ curve $\gamma: I \rightarrow\m$ ($I$ real interval) is
called timelike, lightlike,  spacelike or causal when so it is
$\dot \gamma(s)$ for all $s\in I$. For causal curves, this
definition is extended to include piecewise $C^1$ curves: in this
case, the two limit tangent vectors on the breaks must belong to
the same causal cone. Accordingly, causal curves are called either
{\em future} or {\em past} directed depending on the cone of
$\dot\gamma(s)$.

A smooth curve $\gamma: I \rightarrow\m$ is a {\em geodesic} if
its acceleration vanishes, i.e.,
\begin{equation}
\label{geo} \nabla^L_s\dot \gamma(s) = 0\qquad \mbox{for all $s
\in I$,}
\end{equation}
where $\nabla^L_s$ denotes the covariant derivative along $\gamma$
induced by the Levi--Civita connection of metric $\g_L$. In this
case, the product
\[
E_\gamma \equiv \langle\dot \gamma(s),\dot \gamma(s)\rangle_L
\quad
       \mbox{for all}\; s \in I
\]
is necessarily constant. The spacetime $(\m,\g_L)$ is {\em
geodesically connected} if, given any two points $p, q \in \m$,
there exists a geodesic $z^*: [0,1] \rightarrow\m$ such that
$z^*(0) =p$ and $z^*(1) = q$. This property is equivalent to the
existence of a critical point of the {\em action functional}
defined in (\ref{funzionaleazione}) in the set of all the $C^1$
curves $z: [0,1] \rightarrow\m$ such that $z(0) =p$ and $z(1) = q$
(and also in the extended domain of $H^1$ curves below).

A vector field $K$ in $(\m,\g_L)$ is said {\em complete} if its
integral curves are defined on the whole real line. On the other
hand, $K$ is said {\em Killing} if the Lie derivative of the
metric tensor $\g_L$ with respect to $K$ vanishes everywhere, or,
equivalently, if the stages of all its local flows are isometries
(i.e., $\g_L$ is invariant by its flow).

A well--known characterization of Killing vector fields is the
following: $K$ is a Killing vector field if and only if for each
pair $Y$, $W$ of vector fields, it is
\[
\langle\nabla^L_Y K,W\rangle_L = -\ \langle\nabla^L_W K,Y\rangle_L
.
\]
Thus, if $z:I \to \m$ is a $C^1$ curve and $K$ is a Killing vector
field it is
\[
\langle\dot z,\nabla^L_s K(z)\rangle_L \equiv 0
\]
on $I$. If $z$ is only absolutely continuous, this holds almost
everywhere in $I$. In particular, if $z$ is a geodesic this
property implies the existence of a constant $C_z \in \R$ such
that
\begin{equation}
\label{constraint} \langle\dot z,K(z)\rangle_L \equiv C_z\quad
\hbox{for all $s \in I$.}
\end{equation}

A spacetime is called {\em stationary} if it admits a timelike
Killing vector field $K$. Locally, any stationary spacetime looks
like a {\em standard stationary} one, i.e., the  spacetime in
(\ref{met}). For these spacetimes, without loss of generality it
can be assumed $K=\partial_t$ so to define the (future)
time--orientation. If, in addition, $K$ is also irrotational
(i.e., its orthogonal distribution $K^\perp$ is involutive), the
stationary spacetime is called {\em static}; in this case, locally
it looks like a {\em standard static} one (i.e., its spacetime
metric is the product one obtained in (\ref{met}) with $\delta
\equiv 0$). Notice that a static spacetime may be standard
stationary but not standard static (see Remark \ref{r0} below).
\smallskip

\noindent Given $p$, $q \in \m$ the {\em causality relation} $p <
q$ (respectively {\em chronological relation} $p\ll q$) means that
there exists a future--directed causal (respectively timelike)
curve from $p$ to $q$. Denote by $p \le q$ indistinctly either
$p<q$ or $p=q$. Then, for each $p\in\m$ the {\sl causal future}
$J^+(p)$ and the {\sl causal past} $J^-(p)$ are defined as
\[
J^+(p) =\{q\in\m: p\le q\}\qquad\hbox{and}\qquad J^-(p) =\{q\in\m:
q\le p\}.
\]
Taking into account these relations, the space of piecewise $C^1$
causal curves can be extended in a way appropriate for convergence
of curves (cf. \cite[pp. 442]{Ge} or also \cite[pp. 54]{BEE}):
\begin{definition}\label{dcont}
{\em A (non--necessarily smooth) {\em future-directed causal
curve} $\gamma: I\rightarrow \m $ is a (continuous) curve which,
for each convex\footnote{i.e., $U$ is a (starshaped) normal
neighbourhood of all its points (see \cite[pp. 129]{ON}).}
neighbourhood $U$, satisfies that, given $t,t' \in I, t<t'$ with
$\gamma([t,t']) \subset U$, necessarily $\gamma(t) <_U
\gamma(t')$, where $<_U$ is the causal relation in $U$, regarded
as a spacetime (i.e., $\gamma(t)$ and $\gamma(t')$ can be joined
by a future--directed $C^1$--causal curve contained entirely in
$U$). }\end{definition}

\begin{remark} \label{lipschitz}
{\em Causal curves, even if non--necessarily smooth, must be at
least locally Lipschitzian and, thus, a.e. differentiable  with
finite integral of their length (see \cite[Remark 2.26]{Pe}).
Notice that a continuous curve, a.e. differentiable, with timelike
gradient (in the same time--orientation at each differentiable
point) and finite integral of its length, is not necessarily a
causal curve. A counterexample in Lorentz--Minkowski spacetime
$\LL^2$ can be constructed as follows. Consider a ``devil's
staircase'' type function $t\in[0,1]\mapsto x(t)\in \R$ (the
typical example is Cantor's function) which is continuous, with 0
derivative a.e., and connects $x(0)=0$, $x(1)=2$. Now, the curve
in natural coordinates of $\LL^2$, $\gamma(t)=(x(t),t)$, satisfies
all the required properties, but it connects
the non--causally related points $(0,0), (2,1)$.\\
Recall also that causal curves are absolutely continuous and,
thus, they lie in the spaces of $H^1$ type defined below.}
\end{remark}

There are some equivalent definitions on what means to be globally
hyperbolic for a spacetime:
\begin{itemize}
\item[(1)] The spacetime is strongly causal (i.e.,  no
``almost--closed'' causal curves exist) and $J^+(p)\cap J^-(q)$ is
compact for any $p,q\in \m$. Even more, it is worth pointing out
that the assumption of being strongly causal can be weakened in
only causal (absence of closed causal curves, see \cite{MS} for
detailed explanations).

\item[(2)] The space of causal curves joining any two fixed points
$p, q\in \m$ (defined from $[0,1]$ to $\m$, but identified up to a
strictly increasing monotonic reparametrization) is compact. The
definition of the topology in such space of causal curves is
somewhat subtle (see \cite{Le,Pe}). Essentially, a priori we will
exclude the existence of closed causal curves (otherwise,
parametrizing one such a curve by giving more and more rounds, a
sequence of non--equivalent causal curves would be obtained, and
the compactness of the space of causal curves would be violated)
and, then, the $C^0$ topology of curves is used. \item[(3)] There
exists a Cauchy hypersurface, that is, a subset which is crossed
exactly once by any inextendible timelike curve.
\end{itemize}

A Cauchy hypersurface is necessarily a closed subset of $\m$ and
an embedded topological hypersurface. A long--standing folk
question has been if any globally hyperbolic spacetime must also
admit a smooth Cauchy hypersurface which is spacelike (at all its
points). Recently, this question has been answered affirmatively
in \cite{BeS1} and, thus, we can take  as a characterization of
global hyperbolicity the existence of a {\em (smooth) spacelike
Cauchy hypersurface} $\es \subset \m$. This characterization has
the following remarkable consequence for stationary spacetimes:
\begin{theorem} \label{t0}
A globally hyperbolic stationary spacetime is a standard
stationary one, if some of its timelike Killing vector fields $K$
is complete.
\end{theorem}
\dimo Let $\es$ be a spacelike Cauchy hypersurface, and consider
the map
\[
\Psi: (x,t) \in \es\times\R \mapsto\Phi_{t}(x) \in \m,
\]
where $\Phi$ is the flow of the complete vector field $K$. As each
point of $\m$ is crossed by one integral curve of $K$, which
crosses $\es$ at exactly one point, $\Psi$ is a diffeomorphism. As
$K$ is Killing, the pull--back metric $\Psi^* \g_L$ is independent
of $t$ and, thus, it makes $\es\times\R$ be a standard stationary
spacetime. \cvd

\begin{remark} \label{r0} {\em
(1) If $K$ is also irrotational, Theorem \ref{t0} does not yield
the standard static splitting, as an integral manifold of
$K^\perp$ may be non--Cauchy. A counterexample would be
\[
\m = S^1\times \R,\qquad \g_L=d\theta^2 + 2 dt d\theta -dt^2,
\]
where $(S^1, d\theta^2)$ is the standard unit
circumference, and $K=\partial_t$.

(2) Function $t$ on $\m$ obtained from $\Psi^{-1}$ is a {\em
Cauchy temporal function}, that is, the levels $t=$ constant are
Cauchy hypersurfaces, and $t$ is smooth with a past--directed
timelike gradient (in particular, $t$ is a time function, i.e. a
continuous function which increases on any future--directed causal
curve). As proved in \cite{BeS2}, when such a temporal function
exists the spacetime admits a global orthogonal splitting as in
(\ref{met}) with $\delta \equiv 0$ {\em but with $\beta$ and $\g$
depending on $t$}. This splitting is obtained by flowing through
the integral curves of $\nabla t$ (which, of course, are not equal
to the integral curves of $K\equiv \partial_t$ in general) and,
thus, it has a different nature from the splitting in Theorem
\ref{t0}.

(3) Under the hypotheses of Theorem \ref{globhyp}, $\es$ can be
chosen complete and, thus, so it is $\g$ in the splitting
(\ref{met}).}
\end{remark}
Finally, the following well--known property of globally hyperbolic
spacetimes is stated for reference below (see, for example,
\cite[Proposition 6.6.6]{HE}):

\begin{proposition} \label{jcompact}
If $\m$ admits a Cauchy hypersurface $\es$ then $J^-(p)\cap \es$
is compact, for any $p\in \m$.
\end{proposition}


\section{Abstract intrinsic functional framework}\label{sec3}

Throughout this section we will assume that $(\m,\g_L)$ is a
finite dimensional stationary spacetime with Killing vector field
$K$. Next, geodesic connectedness of $(\m,\g_L)$ will be studied
by using an intrinsic approach and, so, the framework introduced
by Giannoni and Piccione in \cite{GP} is revisited.

In order to define notions as uniform convergence of curves or
$H^1$ spaces, fix any auxiliary Riemannian metric $\g_R$ on $\m$.
This metric can be chosen by leaving $\g_L$ unaltered on the
orthogonal bundle of $K$, and reversing the sign on $K$,
explicitly:
\begin{equation}\label{Riemann}
\langle\zeta,\zeta'\rangle_R = \langle\zeta,\zeta'\rangle_L\ -\ 2\
\frac{\langle\zeta,K(z)\rangle_L\ \langle\zeta',K(z)\rangle_L}
{\langle K(z),K(z)\rangle_L}
\end{equation}
for all $\zeta$, $\zeta'\in T_z\m$, $z\in\m$ (this is the
canonical choice in \cite{GP}).  But recall that, on a standard
stationary spacetime, i.e., when $\m$ is equipped with metric
(\ref{met}), metric (\ref{Riemann}) does not agree with $\g$ on
$\m_0$, and may be incomplete on $\m$. Nevertheless, the results
on this section are independent of the particular choice of
$\g_R$.

As noticed in the previous section, the conservation law
(\ref{constraint}) is a natural constraint for geodesics in
stationary spacetimes. Therefore, it results natural to look for
critical points of action functional $f$ in
(\ref{funzionaleazione}) defined on the set of curves
\[
\begin{array}{rl}
C^1_K(p,q) = \{z\in C^1([0,1],\m) :& z(0) = p,\, z(1) = q,\, \hbox{and $C_z\in \R$ exists}\\
  &\hbox{such that $\langle\dot z,K(z)\rangle_L
\equiv C_z$}\}.  \end{array}
\]
As a first variational principle, we have (see \cite[pp.
2]{GP}\footnote{In \cite[pp. 2]{GP} this results is actually
stated for $C^{1}$ curves (not necessarily in $C^1_K(p,q)$). In
order to pass to $C^1_K(p,q)$ a similar argument to that of
\cite[Proof of Theorem 3.3]{GP} is necessary.}):
\begin{theorem} \label{tcpq}
If $z \in C^1_K(p,q)$ is a critical point of $f$ restricted to
$C^1_K(p,q)$ then $z$ is a geodesic connecting $p$ and $q$.
\end{theorem}

Even if functional $f$ is defined in $C^1_K(p,q)$, it cannot be
managed only in this space, as this space is ``too small'' for
problems of convergence. So, the ``natural'' setting of this
variational problem is a suitable submanifold of the space of
$H^1$ curves from $[0,1]$ to $\m$, named $H^1([0,1],\m)$.

Thus, we define the infinite dimensional manifold
\[\begin{array}{rl}
\Omega^1(p,q) = \{z:[0,1] \to\m:& \hbox{$z$ is absolutely continuous and such that} \\
&\displaystyle z(0) = p,\, z(1) = q,\, \int_0^1\langle\dot z,\dot
z\rangle_R ds < +\infty\},
\end{array}
\]
whose tangent space in each $z \in \Omega^1(p,q)$ can be
identified with
\[\begin{array}{rl}
T_z\Omega^1(p,q) = &\{\zeta:[0,1] \to T\m:\
\zeta(s) \in T_{z(s)}\m\, \hbox{for all $s\in [0,1]$, $\zeta$ is}\\
&\qquad \hbox{absolutely continuous and}\, \zeta(0) = 0=
\zeta(1),\, \|\zeta\|_* < +\infty\},
\end{array}
\]
being its Hilbert norm
\[
\|\zeta\|_*^2 =
\int_0^1\langle\nabla_s^R\zeta,\nabla_s^R\zeta\rangle_R ds,
\]
where $\nabla_s^R$ denotes the covariant derivative along $z$
relative to metric tensor $\g_R$ (nevertheless, we are not
interested in its concrete value, which depends on the chosen
$\g_R$, but only in the finiteness of the norm). Recall that
functional $f$ in (\ref{funzionaleazione}) is well defined and
finite on all $\Omega^1(p,q)$ (for example, notice that, for the
choice (\ref{Riemann}) of $\g_R$, $\langle\zeta,\zeta\rangle_R \ge
|\langle\zeta,\zeta\rangle_L|$). Even more, $f$ is smooth with
differential given by
\[
f'(z)[\zeta] = \int_0^1 \langle\dot z,\nabla_s^L\zeta\rangle_L ds
\]
for all $\zeta \in T_z\Omega^1(p,q)$, $z\in \Omega^1(p,q)$.
Standard calculations allow one to prove that the critical points
in $\Omega^1(p,q)$ are smooth curves which satisfy geodesic
equation (\ref{geo}).

Analogously, the set $C^1_K(p,q)$ can be extended to  a new subset
of $\Omega^1(p,q)$ defined as
\begin{equation}\label{enne}
\begin{array}{rl}
\NN = \{z\in \Omega^1(p,q) : & C_z \in \R\ \hbox{exists such that}\\
& \langle\dot z,K(z)\rangle_L= C_z\ \hbox{a.e. on
$[0,1]$}\}.\end{array}
\end{equation}
In fact, standard arguments in Sobolev spaces imply that the
closure of $C^1_K(p,q)$ is contained in $\NN$, furthermore $\NN$
is a $C^2$--submanifold of $\Omega^1(p,q)$ (see \cite[Proposition
3.1]{GP}) whose tangent space in each point $z \in \NN$ is
\[\begin{array}{rl}
T_z\NN = \{\zeta\in T_z\Omega^1(p,q) :&
\langle\nabla_s\zeta,K(z)\rangle_L + \langle\dot z,\nabla_\zeta K(z)\rangle_L\\
&\hbox{is constant a.e. on $[0,1]$}\}\end{array}
\]
(see \cite[Corollary 3.2]{GP}). For simplicity, denote the
restriction of action functional $f$ on $\NN$ still with $f$.
Thus, the following variational principle can be stated (for a
complete proof, see \cite[Theorem 3.3]{GP}).
\begin{theorem}
A curve $z \in \Omega^1(p,q)$ is a geodesic in $\m$ if and only if
$z \in \NN$ and $z$ is a critical point of functional $f$ in
$\NN$.
\end{theorem}

\begin{remark} \label{lipschitz2}
{\rm An essential step in the proof of Theorem \ref{globhyp} will
be to construct causal curves from certain curves which connect
two fixed points. If these curves belong to a general
$H^1$--space, then we should extend the notion of $C^1$ causal
curve to $H^1$ ones. When one makes this extension some subtleties
appear (recall Remark \ref{lipschitz}), and the space of causal
curves is reobtained (see Appendix B). Nevertheless, we will skip
such technicalities by using sequences of curves in $C^1_K(p,q)$
and their limits in $\NN \subset \Omega^1(p,q)$.}
\end{remark}

Now, let us introduce the following definition which, essentially,
translates classical ``condition $(C)$ of Palais--Smale''  to our
ambient.

\begin{definition} Fixed $c\in \R$ the set $\Omega^1_K(p,q)$ is
$c$--precompact for $f$ if every sequence $(z_n)_n \subset
\Omega^1_K(p,q)$ with $f(z_n) \le c$ has a subsequence which
converges weakly in $\Omega^1(p,q)$ (hence, uniformly in $\m$).
Furthermore, the restriction of $f$ to $\Omega^1_K(p,q)$ is
pseudo--coercive if $\Omega^1_K(p,q)$ is $c$--precompact for all
$c \ge \inf f(\Omega^1_K(p,q))$.
\end{definition}
The geodesic connectivity between each $p$ and $q$ will be a
consequence of the following theorem (see \cite[Theorem 1.2]{GP}).
\begin{theorem}\label{intrinsictheo}
If $C^1_K(p,q)$ is not empty and there exists $c > \inf
f(C^1_K(p,q))$ such that $C^1_K(p,q)$ is $c$--precompact then
there exists at least one geodesic joining $p$ to $q$ in $\m$.
\end{theorem}
The assumption $C^1_K(p,q)$ non--empty must be imposed, because
even if the stationary spacetime is globally hyperbolic it may not
hold (see Appendix A.3 (a) for an explicit counterexample).
Nevertheless, the possibility of $C^1_K(p,q)=\emptyset$ can be
ruled out if $K$ is complete (compare with \cite[Lemma 5.7]{GP}).
\begin{proposition} \label{noempty} Under the hypotheses of Theorem \ref{t0}, for each
$p$, $q \in \m$, it is $C^1_K(p,q) \ne\emptyset$.
\end{proposition}
\dimo Let $p=(x_p,t_p)$, $q = (x_q,t_q)\in \m$ be fixed, and
consider splitting (\ref{met}) ensured by Theorem \ref{t0}. As
$\m$ and, thus, $\es$, is connected, a smooth curve $x :[0,1] \to
\es$ exists joining $x_p$ to $x_q$. Now, compute $t: [0,1] \to \R$
by imposing $t(0)=t_p$ and $\dot t =\ \frac{\langle\delta(x),\dot
x\rangle - C}{\beta(x)}$ (i.e., $\langle (\dot x, \dot t),
\partial_t \rangle_L \equiv C$), where constant $C$ is chosen so to
make $\int_0^1\dot t ds = t_q-t_p$. \cvd
\bigskip


\section{The non--canonical global splitting}\label{sec4}

From now on, suppose that $\m$ has a complete timelike Killing
vector field $K$ and is globally hyperbolic with a complete
spacelike Cau\-chy hypersurface $\es$. By Theorem \ref{t0}, we can
consider that the spacetime is the product $\es \times \R$, with
the metric (\ref{met}) for a certain vector field $\delta$ on
$\es$ and the identifications
\[
K(z) \equiv (0,1) \in T_x\es \times \R\quad \hbox{for all $z = (x,t) \in \m$,}
\]
\[
\langle K(z),K(z)\rangle_L =- \beta(x) \quad \hbox{for all}\quad  z = (x,t)
\in \m.
\]
Nevertheless, recall that neither $K$ nor $\es$ are unique. Thus,
this global splitting is not canonically associated to a spacetime
under hypotheses of Theorem \ref{globhyp}. Anyway, the results
will be independent of the chosen $K, \es$.

For any absolutely continuous curve $z = (x,t) : [0,1] \to \m$, it
is
\begin{equation}\label{kappa}
\langle \dot z(s),K(z(s))\rangle_L = \langle\delta(x(s)),\dot
x(s)\rangle - \beta(x(s))\dot t(s) \quad \hbox{for a.e. $s \in
[0,1]$.}
\end{equation}
Fixed $p = (x_p,t_p)$, $q = (x_q,t_q) \in \m$, it is
\[
\Omega^1(p,q) \equiv \Omega^1(x_p,x_q;\es) \times W(t_p,t_q),
\]
where
\begin{eqnarray*}
&&\begin{array}{rl} \Omega^1(x_p,x_q;\es)
= &\{x:[0,1] \to\es: \hbox{$x$ is absolutely continuous and}\\
&\qquad\qquad\displaystyle x(0) = x_p,\, x(1) = x_q,\,
\int_0^1\langle\dot x,\dot x\rangle ds < +\infty\},
\end{array}\\
&&W(t_p,t_q) = \{t\in H^1([0,1],\R): t(0) = t_p,\, t(1)=
t_q\}=H^1_0([0,1],\R) + T^*,
\end{eqnarray*}
with $H^1([0,1],\R)$ classical Sobolev space and
\[
H_0^1([0,1],\R) =\{t\in H^1([0,1],\R): t(0) = t(1) = 0\},
\]
\[
T^* : s\in [0,1] \longmapsto t_p + s \Delta_t \in\R ,\quad
\Delta_t = t_q - t_p.
\]
Whence, $W(t_p,t_q)$ is a closed affine submanifold of
$H^1([0,1],\R)$ with tangent space 
$T_tW(t_p,t_q) = H_0^1([0,1],\R)$ for all $t \in W(t_p,t_q)$. Moreover, it is
\[
\begin{array}{rl}
T_x\Omega^1(x_p,x_q;\es) = &\{\xi :[0,1]\to T_x\es :
\hbox{$\xi$ is absolutely continuous and}\\
&\qquad\qquad\displaystyle \xi(0) = \xi(1) = 0,\, \int_0^1\langle
D_s \xi,D_s \xi\rangle ds < +\infty\}
\end{array}
\]
for all $x\in\Omega^1(x_p,x_q;\es)$, where $D_s$ denotes the
covariant derivative along $x$ induced by the
Levi--Civita connection of metric $\langle \cdot,\cdot\rangle$.\\
Thus, taken any curve $z = (x,t)\in \Omega^1(p,q)$ it is
\[
T_z \Omega^1(p,q) \equiv T_x\Omega^1(x_p,x_q;\es)\times
H_0^1([0,1],\R)
\]
and $\Omega^1(p,q)$ can be equipped with the Riemannian structure
\[
\langle\zeta,\zeta\rangle_H = \langle(\xi,\tau),(\xi,\tau)\rangle_H =\int_0^1\langle D_s \xi,D_s
\xi\rangle\ ds\ +\ \int_0^1\dot\tau^2 ds
\]
for any $z=(x,t)\in \Omega^1(p,q)$ and $\zeta=(\xi,\tau)\in T_z
\Omega^1(p,q)$.

By Nash Embedding Theorem, Riemannian hypersurface $\es$ can be
assumed as a submanifold of a suitable Euclidean space $\R^N$ with
$\g$ restriction to $\es$ of its Euclidean metric and
$d(\cdot,\cdot)$ the corresponding distance. Furthermore,
$\Omega^1(x_p,x_q;\es)$ is a submanifold of classical Sobolev
space $H^1([0,1],\R^N)$ and is complete because $\es$ is complete.

Clearly, (\ref{enne}) and (\ref{kappa}) imply that $z=(x,t) \in
\NN$ if and only if $x \in \Omega^1(x_p,x_q;\es)$, $t \in
W(t_p,t_q)$ and a constant $C_z \in \R$ exists such that
\begin{equation}\label{time0}
\langle\delta(x),\dot x\rangle - \beta(x) \dot
t=C_z\qquad\hbox{a.e. on $[0,1]$.}
\end{equation}
Hence, it is
\begin{equation}\label{time1}
\dot t =\ \frac{\langle\delta(x),\dot x\rangle -
C_z}{\beta(x)}\qquad\hbox{a.e. on $[0,1]$,}
\end{equation}
which implies
\begin{equation}\label{time2}
C_z =\ \left(\int_0^1\frac{\langle\delta(x),\dot
x\rangle}{\beta(x)}\ ds\ -\ \Delta_t\right)\
\left(\int_0^1\frac{ds}{\beta(x)}\right)^{-1}.
\end{equation}
Moreover, by (\ref{time1}) it follows
\[
\int_0^1\beta(x)\dot t^2 ds =\ \int_0^1\frac{\langle\delta(x),\dot
x\rangle^2}{\beta(x)}\ ds -\ 2 C_z
\int_0^1\frac{\langle\delta(x),\dot x\rangle}{\beta(x)}\ ds +
C_z^2 \int_0^1\frac{ds}{\beta(x)},
\]
and thus, the restricted action functional $f$ becomes
\begin{eqnarray}
f(z) &=& \frac 12 \int_0^1 \langle\dot x,\dot x\rangle ds\
 +\ C_z \Delta_t  +\ \frac 12 \int_0^1 \beta(x)\dot t^2\ ds\nonumber\\
&=& \frac 12 \int_0^1 \langle\dot x,\dot x\rangle ds\ +\ \frac 12
\int_0^1\frac{\langle\delta(x),\dot x\rangle^2}{\beta(x)}\ ds\ -\
\frac{C_z^2}2 \int_0^1\frac{ds}{\beta(x)}.\label{nuovo}
\end{eqnarray}
In conclusion, now we can state the following variational
principle introduced by Giannoni and Masiello in \cite{giama1} for
standard stationary spacetimes (see also \cite[Theorem
3.3.2]{Ma}):
\begin{proposition}
\label{newazione} A curve $z^* = (x^*,t^*) \in \Omega^1(p,q)$ is a
critical point of action functional $f$ in $\Omega^1(p,q)$ if and
only if $x^*$ is a critical point of functional
\[
\J: \Omega^1(x_p,x_q;\es) \to \R
\]
defined in (\ref{newfunction}) and $t^* = \Psi(x^*)$ with
$\Psi:\Omega^1(x_p,x_q;\es) \to W(t_p,t_q)$ such that
\begin{equation}\label{time5}
\begin{array}{l}
\Psi(x)(s) = \displaystyle t_0 + \int_0^s\frac{\langle
\delta(x(\sigma)),\dot x(\sigma)\rangle}{\beta(x(\sigma))}\ d\sigma \\
\qquad\displaystyle -\ \left(\int_0^1 \frac{\langle \delta(x),\dot
x\rangle}{\beta(x)}\ ds - \Delta_t\right)\ \int_0^s
\frac1{\beta(x(\sigma))}\ d\sigma\ \left(\int_0^1
\frac1{\beta(x)}\ ds\right)^{-1}.
\end{array}
\end{equation}
Moreover, it is $f(z^*) = \J(x^*)$.
\end{proposition}

\begin{remark}
{\rm Taken any $z =(x,t) \in \NN$ formulae (\ref{newfunction}),
(\ref{time2}) and (\ref{nuovo}) imply that $f(z) = \J(x)$ while
(\ref{time1}), (\ref{time2}) and (\ref{time5}) imply $t =\Psi(x)$.
In particular, let us point out that in the proof of Proposition
\ref{newazione} the introduced variational constraint is the
kernel of the ``partial derivative'' operator
\[
z \mapsto f'(z)[(\cdot,0)]
\]
which gives exactly the natural constraint (\ref{time0}).}
\end{remark}


\section{Geodesic connectedness} \label{sec5}

The aim of the present section is to prove the following result.
\begin{theorem}
\label{coerciveness} Under the hypotheses of Theorem
\ref{globhyp},  the restriction of $f$ to $C^1_K(p,q)$ is
pseudo--coercive, for any $p$, $q \in \m$.
\end{theorem}
Thus, Theorem \ref{globhyp} follows directly from the intrinsic
result on geodesic connectedness for precompact $C^1_K(p,q)$ (see
Theorem \ref{intrinsictheo}), the  non--emptyness of $C^1_K(p,q)$
(see Proposition \ref{noempty}), and the coercivity ensured in
Theorem \ref{coerciveness}.

The key step for the proof of Theorem \ref{coerciveness} is the
following proposition, which relates technically the global
hyperbolicity of the spacetime with the notion of
pseudo--coercivity. Let us remark that, in what follows, the
hypotheses of Theorem \ref{globhyp} are always assumed, as well as
the standard stationary splitting associated to the complete $K$
and ${\cal S}$, and the notations introduced in Sections
\ref{sec3} and \ref{sec4}. In particular, for any curve $x$ in
${\cal S}$, $\|\dot x\|^2 = \int_0^1\langle\dot x,\dot x\rangle
ds$, where the metric $\langle \cdot, \cdot \rangle$ is just the
induced Riemannian metric on ${\cal S}$.
\begin{proposition} \label{p-s}
Let $(z_n)_n$, $z_n=(x_n,t_n)$, be a sequence of curves in
$C^1_K(p,q)$ such that $f(z_n) (=\J(x_n))$ is upper bounded for
all $n$. Then, $(\|\dot x_n\|)_n$ is bounded, too.
\end{proposition}

Now, let us state some lemmas useful in the proof of Proposition
\ref{p-s} in which, for simplicity, we define 
$C^1(x_p,x_q;\es) = \Omega^1(x_p,x_q;\es) \cap C^1([0,1], \m)$.

Firstly, let us point out that when the $x_n$'s lie in a compact
subset of ${\cal S}$, Proposition \ref{p-s} is just a consequence
of Cauchy--Schwarz inequality.
\begin{lemma}\label{lemma-1}
Let $(x_n)_n \subset C^1(x_p,x_q;\es)$. If a compact subset ${\cal
C}$ of ${\cal S}$ contains all the elements of the sequence
$(x_n)_n$ and $\|\dot x_n\| \to +\infty$, then $\J(x_n) \to
+\infty$.
\end{lemma}
\dimo Consider definition (\ref{newfunction}). By expanding the
squared term and using Cauchy--Schwarz inequality:
\begin{equation}\label{lowbound}
2\J(x_n) \ge \|\dot x_n\|^2 - \Delta_t \left(\Delta_t - 2
\int_0^1\frac{\langle\delta(x_n),\dot x_n\rangle}{\beta(x_n)}\
ds\right) \left( \int_0^1\frac{ds}{\beta(x_n)}\right)^{-1}
\end{equation}
for each $n\in\N$. Now, from the compactness of ${\cal C}$, there
exist some strictly positive constants $N_1$, $N_2$, $\nu$ such that,
\begin{equation}\label{bounds2}
\nu \le \beta(x_n(s))\le N_1 ,\,
\langle\delta(x_n(s)),\delta(x_n(s))\rangle \le N_2 , \, \hbox{for
all $s \in [0,1]$, $n\in\N$}.
\end{equation}
Thus, by applying again Cauchy--Schwarz, it is
\[
2\J(x_n) \ge \|\dot x_n\|^2 - N \|\dot x_n\| -N'
\]
for some $N, N' >0$ independent of $n$, and the result follows.
\cvd
\bigskip

On the contrary, when no such compact ${\cal C}$ exists,
Proposition \ref{p-s} will be proved in two steps, Lemmas
\ref{lemma0} and \ref{lemma2}. But first notice that, given the
spacelike parts $(x_n)_n$, a sequence of future--directed
lightlike curves from $p=(x_p,t_p)$ to the integral curve of $K$
through $q=(x_q,t_q)$ can be constructed. More precisely:
\begin{lemma}\label{lemma1}
Fixed any $x\in C^1(x_p,x_q;\es)$ ($x$ non--constant if $x_p=x_q$)
there exists a unique lightlike curve $\gamma^l =(x^l,t^l): [0,1]
\to \m$ joining $(x_p,t_p)$ to $\{x_q\}\times \R$ in a time $
T(x)= t^l(1)-t^l(0)>0$ such that $x^l=x$. Moreover, $T(x)$
satisfies:
\begin{equation}
\label{time3} T(x) = \int_0^1\frac{\langle\delta(x),\dot
x\rangle}{\beta(x)}\ ds\ +
\int_0^1\frac{\sqrt{\langle\delta(x),\dot x\rangle^2 + \langle\dot
x,\dot x\rangle \beta(x)}}{\beta(x)}\ ds.
\end{equation}
\end{lemma}
\dimo Fixing $x\in C^1(x_p,x_q;\es)$, the $t^l$ part of the curve
$\gamma^l =(x,t^l): [0,1] \to \m$ is characterized by the
equalities $\langle\dot\gamma^l,\dot\gamma^l\rangle_L = 0$,
$t^l(0) = t_p$ and the inequality $\dot t^l > 0$ a.e. on [0,1].
From (\ref{met}) it follows
\[
\dot t^l = \frac{\langle\delta(x),\dot x\rangle +
\sqrt{\langle\delta(x),\dot x\rangle^2 + \langle\dot x,\dot
x\rangle \beta(x)}}{\beta(x)}\quad \hbox{a.e. on [0,1]},
\]
which implies directly (\ref{time3}).\cvd
\bigskip

The global hyperbolicity of $\m$ becomes crucial for the first
conclusion of the following result. The second one is just a
simple consequence of the completeness of ${\cal S}$, but this
property also turns out to be essential.
\begin{lemma}\label{lemma0}
Let $(x_n)_n \subset C^1(x_p,x_q;\es)$ and, for each $n\in\N$,
denote $T_n=T(x_n)$. If no compact subset of ${\cal S}$ contains
all the elements of sequence $(x_n)_n$, then, up to a subsequence,
it is:
\begin{itemize}
\item[{\sl (i)}] $T_n\rightarrow +\infty$, and \item[{\sl(ii)}]
$\|\dot x_n\| \to +\infty$.
\end{itemize}
\end{lemma}
\dimo {\sl (i)} Arguing by contradiction, let $T^+$ be  an upper
bound for all $T_n$, and put $p=(x_p,t_p)$, $q^{+}=(x_{q},T^+)$.
The lightlike curves $\gamma^l_n=(x_n,t_n^l)$ obtained from Lemma
\ref{lemma1}, can be prolongued with the integral curve of $K$
from $(x_{q},T_n)$ to $(x_{q},T^+)$; so, a piecewise
smooth future--directed causal curve from $p$ to $q^{+}$ is
obtained and all the $\gamma_n^l$ lie in $J^-(q^{+})$. By
Proposition \ref{jcompact}, $J^-(q^{+})\cap {\cal S}$ is compact,
but this is a contradiction because this subset contains
all the $x_n$'s.\\
{\sl (ii)} As $\es$ is complete, no bounded subset can contain all
the $x_n$'s. So, there is a sequence of points $(x_n(s_n))_n$ at
arbitrary large distance from $x_p$, and the result follows. \cvd

\begin{lemma}\label{lemma2}
Fixed any sequence $(x_n)_n \subset C^1(x_p,x_q;\es)$ such that
\begin{equation}\label{limit}
\|\dot x_n\| \to +\infty\qquad \hbox{and}\qquad T_n \to +\infty
\end{equation}
then
\begin{equation}\label{limit2}
\J(x_n) \to +\infty.
\end{equation}
\end{lemma}
\dimo Let $(x_n)_n \subset C^1(x_p,x_q;\es)$ be a sequence such
that (\ref{limit}) holds and, for simplicity, let us assume
$\Delta_t > 0$. Taking into account inequality (\ref{lowbound}),
the desired limit (\ref{limit2}) follows from (\ref{limit}), if a
constant $k> 0$ exists such that
\[
\left(\Delta_t - 2 \int_0^1\frac{\langle\delta(x_n),\dot
x_n\rangle}{\beta(x_n)}\ ds\right) \left(
\int_0^1\frac{ds}{\beta(x_n)}\right)^{-1} \le k \quad \hbox{for
all $n \in \N$.}
\]
So, assume that, up to subsequences, it is
\begin{equation}\label{limit3}
\left(\Delta_t - 2 \int_0^1\frac{\langle\delta(x_n),\dot
x_n\rangle}{\beta(x_n)}\ ds\right) \left(
\int_0^1\frac{ds}{\beta(x_n)}\right)^{-1} \ \longrightarrow\
+\infty\quad \hbox{as $n \to +\infty$.}
\end{equation}
On the other hand, by Cauchy--Schwarz inequality and definition
(\ref{time3}), for each $n \in \N$ it is
\begin{equation}\label{5.4b} T_n \le \tilde T_n
\end{equation}
with
\begin{eqnarray*}
\tilde T_n &=& \int_0^1\frac{\langle\delta(x_n),\dot x_n\rangle}{\beta(x_n)}\ ds\\
&&+ \sqrt{\left( \int_0^1\frac{\langle\delta(x_n),\dot
x_n\rangle^2}{\beta(x_n)}\ ds + \|\dot x_n\|^2\right)
\int_0^1\frac{ds}{\beta(x_n)}}\
\end{eqnarray*}
($\tilde T_n$ is the arrival time $\Delta_t \geq 0$ in the
expression of $\J$ that we must choose in order to obtain
$\J(x_n)=0$; this arrival time is also useful in the context of
Fermat principle as stated in \cite{FGM95}). Thus, it is
\[
\begin{array}{l}
\displaystyle\int_0^1\frac{\langle\delta(x_n),\dot x_n\rangle^2}
{\beta(x_n)}\ ds +
\|\dot x_n\|^2\\
\qquad \displaystyle =\ \left(\tilde T_n -
\int_0^1\frac{\langle\delta(x_n),\dot x_n\rangle}{\beta(x_n)}\
ds\right)^2 \left(\int_0^1\frac{ds}{\beta(x_n)}\right)^{-1}\qquad
\hbox{with $\tilde T_n \to +\infty$.}\end{array}
\]
Hence, from (\ref{newfunction}) it follows
\begin{eqnarray*}
2\J(x_n) &=& \left(\tilde T_n -
\int_0^1\frac{\langle\delta(x_n),\dot x_n\rangle}{\beta(x_n)}\
ds\right)^2
\left(\int_0^1\frac{ds}{\beta(x_n)}\right)^{-1}\\
&&-\ \left(\int_0^1\frac{\langle\delta(x_n),\dot
x_n\rangle}{\beta(x_n)}\ ds\ -\ \Delta_t\right)^2
\left(\int_0^1\frac{ds}{\beta(x_n)}\right)^{-1}\\
&=& \left(\tilde T_n^2 - \Delta_t^2 - 2 (\tilde T_n - \Delta_t)
\int_0^1\frac{\langle\delta(x_n),\dot x_n\rangle}{\beta(x_n)}\
ds\right)\
\left(\int_0^1\frac{ds}{\beta(x_n)}\right)^{-1}\\
&=& (\tilde T_n - \Delta_t)\ \left(\tilde T_n + \Delta_t - 2
\int_0^1\frac{\langle\delta(x_n),\dot x_n\rangle}{\beta(x_n)}\
ds\right)\ \left(\int_0^1\frac{ds}{\beta(x_n)}\right)^{-1}.
\end{eqnarray*}
In conclusion, (\ref{limit2}) follows from (\ref{limit}),
(\ref{limit3}) and (\ref{5.4b}). \cvd
\bigskip

\noindent {\bf Proof of Proposition \ref{p-s}.}\\
The proof follows directly from Lemmas \ref{lemma-1},
\ref{lemma0}, \ref{lemma2}. \cvd
\bigskip

\noindent {\bf Proof of Theorem \ref{coerciveness}.}\\
Let $p, q\in \m$ and $c \ge \inf f(C^1_K(p,q))$ be fixed, and
consider a sequence $(z_n)_n$ in $C^1_K(p,q)$ such that
\[
 f(z_n) \le c\quad \hbox{for all $n \in \N$.}
\]
By the global splitting of $\m$, it is $p=(x_p,t_p)$ and
$q=(x_q,t_q)$, while $z_n = (x_n,t_n)$ is such that $x_n \in
C^1(x_p,x_q;\es)$ and $t_n \in C^1([0,1],\R)\cap W(t_p,t_q)$ with
\[
\langle\delta(x_n),\dot x_n\rangle - \beta(x_n) \dot t_n \equiv
C_n\qquad\hbox{on [0,1].}
\]
By Proposition \ref{p-s},
\begin{equation}\label{bdd}
(\|\dot x_n\|)_n\quad \hbox{has to be bounded,}
\end{equation}
and all the $x_n$'s lie in a bounded  subset of $\es$. Thus,
$(x_n)_n$ is bounded in $H^1([0,1],\R^N)$ and, up to subsequences,
it converges to some $x\in H^1([0,1],\R^N)$ weakly in
$H^1([0,1],\R^N)$ and uniformly in [0,1]. As $\es$ is complete,
then $x \in \Omega^1(x_p,x_q;\es)$, and all the $x_n$'s lie in a
compact subset. Thus, (\ref{bounds2}) holds, and, by (\ref{time2})
and Cauchy--Schwarz inequality, sequence $(C_n)_n$ has to be
bounded. Hence, (\ref{time1}) and, again, Cauchy--Schwarz
inequality imply that $(t_n)_n$ is bounded in $H^1([0,1],\R)$ too,
and thus, $t \in W(t_p,t_q)$ exists so that $t_n \to t$ uniformly
in [0,1] (up to subsequences). \cvd


\section*{Appendix A: Discussion on the hypotheses and some counterexamples}

\textbf{\em $\S$1. Estimates for the global hyperbolicity of a
standard stationary spacetime.} As proven in \cite{Sa1}, the
spacetime (\ref{met}) is globally hyperbolic if $\g$ is complete,
$\beta$ is at most quadratic and $\delta$ at most linear, that is
(\ref{quad}) and (\ref{lineare}) hold. This is a ``rough
estimate'': it is easy to find counterexamples when the exponent
of $d$ in any of the inequalities is increased a bit, but these
inequalities are only sufficient conditions. In fact, in the
standard static case $\delta \equiv 0$, the spacetime is globally
hyperbolic if and only if the conformal metric
\begin{equation}\label{eapp} 
\g^*= \frac{\g}{\beta(x)}
\end{equation}
is complete (see \cite{Sa2} for more details). Notice that it is
not relevant for $\g$ to be  complete or not. In fact, classical
Schwarzschild spacetime is globally hyperbolic with incomplete
$\g$. Nevertheless, when a standard static spacetime is globally
hyperbolic the slices at constant $t$ are Cauchy hypersurfaces.
Moreover, such a spacetime admits a {\em complete} spacelike
Cauchy hypersurface $\es$ if and only if $\g$ is complete (the
non--trivial implication to the right can be proven because the
projection $\es \to \mo$ is a diffeomorphism which increases the
distances).
\smallskip

\noindent \textbf{\em $\S$2. Arbitrariness of  standard stationary
splittings.} As explained in the Introduction, essentially all
the previous results in the literature on geodesic connectedness
of stationary spacetimes rely in the  behaviour of $\g, \beta,
\delta$. Nevertheless, these elements are not canonical for the
spacetime, in the sense that many such splittings are possible
with a very different behavior for them. For example, $\LL^2$ can
be written as $\R^2$ endowed with the metric $dx^2+2 \bar
\delta(x) dxdt - dt^2$ for any function $\bar \delta$. This can be
checked because the spacetime is a flat spaceform, i.e., its Gauss
curvature is 0, it is simply connected and geodesically complete.
The last property follows because, as a consequence of
\cite[Proposition 2.1]{RS}, one has just to prove that the metric
$\g_R$ in (\ref{Riemann}) is complete. But this is straightforward
because, in the natural coordinates, the matrix of this metric has
eigenvalues greater than a positive constant (see also
\cite[Example 2.4]{RS}). Thus, vector field $\delta= \bar \delta
\partial_x$ may not satisfy the (at most) linear condition (\ref{lineare}) above.
Even more, it is also easy to construct incomplete Cauchy
hypersurfaces in $\LL^2$ and, then, by using the natural Killing
vector field $K=\partial_t$, Theorem \ref{t0} yields a stationary
splitting with incomplete $\g$.
\smallskip

\noindent \textbf{\em $\S$3. Accuracy of the hypotheses of Theorem
\ref{globhyp}.} Let us check this with two counterexamples:
\begin{itemize}
\item[(a)] {\em Stationary  + Globally hyperbolic with complete
$\es$ $\not\Rightarrow$ geodesically connected}. Consider the
spacetime obtained by removing in Lorentz--Min\-kow\-ski
$\LL^{n+1}, n\geq 1,$ the causal future of the points with
$x_1=0=t$, in natural coordinates $(x_1, \dots , x_n, t)$.
Clearly, the spacetime admits the hyperplane $t\equiv -1$ as a
complete Cauchy hypersurface, but it is not geodesically
connected. Moreover, $C^1_K(p,q)=\emptyset$ for $p=(-1,0,\dots,
0), q=(1,0,\dots , 0)$. \item[(b)] {\em Stationary with complete
$K$ + Globally hyperbolic $\not\Rightarrow$ geodesically
connected}. Let $(\es,\langle\cdot,\cdot\rangle)$ be a
non--geodesically connected Riemannian mani\-fold. Take $\beta$
such that $\g^*$ in (\ref{eapp}) is complete. Then, the standard
static spacetime $(\es\times\R,\g-\beta(x) dt^{2})$ is globally
hyperbolic and not geodesically connected, because the slice $t=0$
is totally geodesic (the family of static spacetimes given in
\cite[Section 7]{BCFS} also stresses the importance of global
hyperbolicity).
\end{itemize}

\noindent \textbf{\em $\S$4. Accuracy of pseudo--coercivity from a
technical viewpoint.} As shown in Theorem \ref{intrinsictheo},
pseudo--coercivity  of functional $f$ (on $C^1_K(p,q)$ for all
$p$, $q\in \m$) yields a technical natural condition for the
geodesic connectedness of the spacetime. Let us discuss the
relation between this condition and others involved in Theorem
\ref{globhyp} as well as  in \cite{GP}. Along this discussion,
$(\m, \g_L)$ is a stationary spacetime with a given timelike
Killing vector field $K$, and we emphasize that, in the point (c),
the Riemannian metric $\g_R$ on $\m$ will be chosen as the one
associated to $\g_L$ and $K$ by formula (\ref{Riemann}).
\begin{itemize}
\item[(a)] {\em Functional $f$ pseudo--coercive $\Rightarrow$ $\m$
globally hyperbolic} (and, thus, the  spacetime admits a splitting
as standard stationary if $K$ is complete). A proof (valid in the
case $K$ complete) can be seen in \cite[Proposition B.1]{GP}. An
alternative argument based on the definitions of global
hyperbolicity in Section \ref{section2} is the following. Arguing
by contradiction, if the space of causal curves joining two points
$p \leq  q$ is not compact, then there exists a sequence $(z_n)_n$
of such future--directed causal curves with no converging
subsequence. As each $z_n$ can be approximated by piecewise smooth
lightlike curves, we can assume that the $z_n$'s are in fact
lightlike curves. Thus, $f(z_n)=0$ for any reparametrization of
$z_n$. As $\langle\dot z_n,K(z)\rangle_L <0$, we can choose this
reparametrization (and smooth the possible finite number of breaks
of $z_n$) in order to make $z_n$ to belong to
$C^1_K(p,q)$ with bounded $f(z_{n})$. So, a converging subsequence of $(z_n)_n$ has to exist, which yields the contradiction.\\
Even more, if $\m$ is standard static then the pseudo--coercivity
also implies $\g$ complete (notice that this does not hold in the
stationary case, as explained at the end of $\S$2). In fact,
otherwise an incomplete geodesic $x: [0,1) \rightarrow \m_0$
exists. The curve in $\m$, $z=(x,0)$ will be a geodesic too, with
$C_z=0$. Now, consider the sequence of curves $z_n = (x_n,0)$,
where each $x_n$ is a loop obtained by reparametrizing the
restriction $x|_{[0, 1-1/n]}$ in such a way that
$p=x(0)=x_n(0)=x_n(1)$, $x_{n}(1/2)=x(1-1/n)$ for all $n>1$.
Again, $C_{z_n}=0$, and $(z_n)_n$ violates the pseudo--coercivity
of $C_{K}^{1}(p,p)$. \item[(b)] {\em Functional $f$
pseudo--coercive $\not\Rightarrow$ $K$ complete.} A simple
counterexample is any strip $\m_0 \times (a,b)$,
$(a,b)\varsubsetneq\R$ of any standard static spacetime $\m=\m_0
\times\R$ with a complete Cauchy hypersurface (for example, $\m=
\LL^n$). In fact, $f$ is pseudo--coercive for the full $\m$ (from
the proof of Theorem \ref{globhyp}), and the strip is still
pseudo--coercive, as any curve $(x,t)\in C^1_K(p,q)$ has $t$
either non--decreasing or
non--increasing.\\
Nevertheless, recall that this is the unique example in the
standard static case, and it does not yield new interesting
examples in the stationary one. In fact, in the stationary case,
when there exists a curve $(x,t)\in C^1_K(p,q)$ (for example, a
geodesic) such that $t$ admits either a strict maximum or minimum,
then no strip $\m_0 \times (a,b)$ with $-\infty < a < b < +\infty$
is pseudo--coercive. 
\item[(c)] {\em Functional $f$
pseudo--coercive + $K$ complete $\Rightarrow$ $\g_R$ in
(\ref{Riemann}) complete} (in particular, $\g_R$ will be complete
under the hypotheses of Theorem \ref{globhyp}). From the final
discussion in (a), the result is obvious in the standard static
case, because then $\g$ becomes complete and $\g_R$ becomes the
Riemannian warped product $\g_R=\g + \beta(x) dt^2$, which is
complete, too (see \cite[Lemma 7.40]{ON}). For the general case,
choose any incomplete $\g_R$--geodesic $\gamma: I \rightarrow \m$
and consider the map $\psi: I\times \R \rightarrow \m$,
$(s,t)\mapsto \Phi_t(\gamma(s))$, where $\Phi$ denotes the flow of
$K$. Then, $I\times \R$, with the induced metric $\psi^*\g_L$, is
static ($\partial_t$ is a timelike Killing vector field and, in
dimension 2, any such vector field is irrotational). Moreover,
$\psi^*\g_R$ coincides with the Riemannian metric on $I\times \R$
obtained from $\psi^*\g_L$ and $\partial_t$ in (\ref{Riemann}).
Recall also that $f$ becomes pseudo--coercive for this spacetime;
thus, $I\times \R$ becomes globally hyperbolic and (as $I\times
\R$ is 2--dimensional and simply connected) standard static. But,
then, the standard static case is applicable, and the metric
$\psi^*\g_R$ must be complete, a contradiction.
\end{itemize}

It is worth pointing out that the completeness of $\g_R$ implies
the completeness of $K$ (if an integral curve of $K$  escapes any
compact subset, it must have infinite length by the completeness
of $\g_R$ and will be complete as it has constant $\g_R$--speed).
Nevertheless, the completeness of $\g_R$ and global hyperbolicity
are independent: any compact stationary spacetime is a
counterexample for the implication to the right, and Schwarzschild
spacetime is a counterexample for the converse.

Finally, we can wonder if both conditions together, the
completeness of $\g_R$ and the global hyperbolicity of $\m$, would
imply the existence of a complete Cauchy hypersurface (and, thus,
geodesic connectedness, by Theorem \ref{globhyp}). Nevertheless,
this type of questions involves completely different techniques
(see, for example, \cite{BeS3}) and goes beyond the scope of the
present article.

\subsection*{Appendix B: $H^1$ causal curves are causal curves}

The most general notion of causal curve is the one for {\em
continuous}  curves (see Definition \ref{dcont}). Necessarily,
such a causal curve is $H^1$ (see Remark \ref{lipschitz}) with
a.e. causal gradient, in the same time--orientation. Now, our
purpose is to prove the converse, which becomes natural when
properties of completeness of the set of causal curves are
considered (see Remark \ref{lipschitz2}). In the proof, the
absolute continuity of $H^1$ curves will be essential (recall the
Cantor--type counterexample in Remark \ref{lipschitz}). So, our
final result can be stated as follows.

\begin{theorem}\label{tappb}
Let $(\m, \g_{L})$ be a spacetime, and $\gamma: [0,1] \rightarrow
\m$, $I=[0,1]$, a continuous curve. Then, the following items are
equivalent:
\begin{itemize}
\item[$(i)$] $\gamma$ is future--directed  causal (according to
Definition \ref{dcont}); \item[$(ii)$] $\gamma$ is $H^1$ and $\dot
\gamma(s)$ is a future--directed causal vector for a.e. $s\in I$.
\end{itemize}
\end{theorem}
\dimo According to Remark \ref{lipschitz}, we have only to prove
that $(ii)\Rightarrow (i)$. Thus, taken $\gamma$ so that $(ii)$
holds, it is enough to show that, chosen any $s_0\in I$ and any
convex neighborhood $U$ of $z_0=\gamma(s_0)$, there exists $0<
\delta <1$
such that if $0<s_1-s_0<\delta$ then $z_0 <_U z_1=\gamma(s_1)$.\\
Firstly, notice that, in a neighbourhood $V\subseteq U$ of $z_0$,
the spacetime can be written as $S\times (t_-,t_+)$ with
$$
\g_{L}[x,t]= g_t[x]-dt^2,
$$
where $g_t$ is a Riemannian metric on $S$ for any $t\in (t_-,t_+)$
and $\partial_t$ is future--directed (for example, see \cite[pp.
53]{Pe}). Then, we can choose $\delta > 0$ such that
$\gamma([s_0,s_0+\delta]) \subset V$.\\
Put $z_i=(x_i,t_i)$, $i=0,1$. Clearly, it is $z_0<_Uz_1$ if there
exists a $C^1$ curve $y: t\in [t_0,t_1]\mapsto y(t) \in S$ such that
$y(t_i)=x_i$, $i=0,1$, and
$$
g_t(\dot y(t), \dot y(t)) \leq 1 \quad \quad \hbox{for all $t\in
[t_0,t_1]$}
$$
(notice that $t\mapsto (y(t),t)$ would be future--directed causal). \\
To this aim, we need the following lemmas.
\begin{lemma} \label{appbl1}
If there exists a sequence of $C^1$ curves $t\in [t_0,t_1]\mapsto
y_n(t) \in S$ such that $ y_n(t_i)=x_i$, $i=0,1$ and, for some
sequence $\epsilon_n \searrow 0$,
$$
g_t(\dot y_n(t), \dot y_n(t)) \leq 1+\epsilon_n \quad \quad
 \hbox{for all $t\in [t_0,t_1]$,}
$$
then $z_0<_Uz_1$.
\end{lemma}
\dimo The curves $t\mapsto
(y_n(t),t_0+\sqrt{1+\epsilon_n}(t-t_0))$ are future--directed and
causal; thus, $z_0<_U z_{1,n}= (x_1,
t_0+\sqrt{1+\epsilon_n}(t_1-t_0))$. As $U$ is convex, the relation
$\leq_U$ is closed and the result follows passing to the limit as
$n\rightarrow +\infty$. \cvd

\begin{lemma} \label{appbl2}
Writing the $H^1$ curve $\gamma$ as $\gamma(s)=(x(s),t(s))$, if
$s\in[s_0,s_0+\delta]$, then function $s\mapsto t(s)$ is strictly
increasing.
\end{lemma}
\dimo As $\gamma$ is absolutely continuous, then
$$
t(s_3)=t(s_2) +\int_{s_2}^{s_3} \dot t(s) ds\qquad\hbox{for all
$s_2$, $s_3 \in [s_0,s_0+\delta]$.}
$$
But $\dot \gamma (s)$ is a.e. future--directed, so $\dot t(s)>0$
a.e. in $I$; hence, $t(s_2) < t(s_3)$ if $s_2<s_3$. \cvd
\bigskip

By Lemma \ref{appbl2} it follows that the $x$--component of
$\gamma$ can be reparametrized by $t$, and this reparametrized
curve still belongs to $H^{1}$. In fact, it is enough to prove
that $t\mapsto x(t)\equiv x(s(t))$ satisfies Barrow's rule
$x(t)=x(t_0) + \int_{t_0}^t (dx/dt)(\bar t)d\bar t$, and this is
obvious because $dx/dt= (dx/ds) (ds/dt)$ a.e., the change of
variable theorem for a.e. differentiable functions is applicable
(see, for example, \cite[Theorem 6.94]{St}), and function
$s\mapsto x(s)$ does satisfy Barrow's rule. Moreover, as
$\dot\gamma$ is future--directed, it has to be
\begin{equation}\label{*}
g_t\left(\frac{dx}{dt}, \frac{dx}{dt}\right) \leq 1 \quad \quad
\mbox{a.e. in $[t_0,t_1]$.}
\end{equation}
Now, take any sequence $(\tilde{y}_{n})_{n}$ of $C^{1}$ curves in
$S$ with $\tilde{y}_{n}(t_{i})=x_{i}$, $i=0,1$, which approach
$t\mapsto x(t)$ in the $H^{1}$ norm. In particular,
$(\tilde{y}_{n}')_{n}$ also go to $x'$ strongly in $L^{2}$ norm.
In order to check that the sequence of $C^{1}$ curves
$(y_{n})_{n}$ formed by the reparametrizations of $\tilde{y}_{n}$
to constant speed, falls necessarily under the hypotheses of Lemma
\ref{appbl1}, firstly notice that the lengths $L_n =$ length$(\tilde
y_n) =$ length$(y_n)$ satisfy $L_n\rightarrow L=$ length$(x)$. Thus,
$$
\left|\frac{dy_n}{dt} \right| = \frac{L_n}{t_1-t_0} \rightarrow
\frac{L}{t_1-t_0} \leq 1
$$
(the last inequality from (\ref{*})), and Lemma \ref{appbl1}
applies.

 \cvd

\subsection*{Acknowledgments}

The authors would like to acknowledge the International Erwin
Schr\"odinger Institute for Mathematical Physics (ESI) in Vienna,
where this work started, for its kind hospitality during the {\sl
Special Research Semester on ``Geometry of pseudo--Riemannian
Manifolds with Application in Physics''} organized by Dmitri
Alekseevski, Helga Baum and Jerzy Konderak.


\end{document}